\documentclass{birkjour}

\usepackage{amsmath}
\usepackage{hyperref}

\newcommand{\Hil}{\mathcal{H}}

\newtheorem{thm}{Theorem}[section]
\newtheorem{cor}[thm]{Corollary}
\newtheorem{lem}[thm]{Lemma}
\newtheorem{prop}[thm]{Proposition}
\theoremstyle{defn}

\theoremstyle{remark}
\newtheorem{rem}[thm]{Remark}
\newtheorem*{ex}{Example}
\numberwithin{equation}{section}

\title[Dynamical frames through Sz.-Nagy–Foiaş dilations]
{Dynamical frames \\ through Sz.-Nagy–Foiaş dilations}

\author{Sebastian Gwizdek}
\address{AGH University of Krakow,
Faculty of Applied Mathematics,
Al. A. Mickiewicza~30, 30-059 Krakow,
Poland}
\email{sebastian.gwizdek@agh.edu.pl}

\author{Sergiusz Ku\.{z}el}
\address{(corresponding author)\br  
AGH University of Krakow,
Faculty of Applied Mathematics,
Al. A. Mickiewicza~30, 30-059 Krakow,
Poland}
\email{kuzhel@agh.edu.pl}

\author{Witold Majdak}
\address{AGH University of Krakow,
Faculty of Applied Mathematics,
Al. A. Mickiewicza~30, 30-059 Krakow,
Poland}
\email{majdak@agh.edu.pl}

\subjclass{42C15, 47A20, 47A45, 47A48}

\keywords{contraction, dilation, characteristic function, dynamical frames}

\begin{document}

\begin{abstract}
We present a unified operator-theoretic treatment of dynamical frames based on defect operators, Stein equations, and the Sz.-Nagy--Foiaş dilation theory of contractions. An explicit realization of the minimal isometric dilation relates generating vectors to the wandering subspace and provides a geometric interpretation of standard frame constructions. For Parseval dynamical frames, polar decomposition yields a complete factorization of the generating families through the adjoint defect operator. General dynamical frames are described in terms of positive solutions of a Stein equation, leading to an equivalent formulation of the frame index. We also collect characteristic-function criteria for strong stability and joint strong stability. In the diagonal case, the characteristic function expresses the classical separation and Carleson conditions in operator-theoretic terms. The emphasis is on explicit constructions, direct proofs, and a common framework for known~results.
\end{abstract}

\maketitle
   
\section{Introduction}
Over the last decade, the theory of dynamical frames
has attracted increasing attention; see the survey article \cite{Rev} and the references therein. 
One of the principal motivations for studying such frames arises from the theory of dynamical sampling. 
Dynamical sampling, introduced in \cite{first}, 
is motivated by the observation that sampling the solution 
of a dynamical system at discrete locations and successive time instants can provide sufficient information to 
reconstruct the underlying system or its initial state. 
From a mathematical point of view, in the discrete-time
setting, this reconstruction problem can be reformulated 
as the question of whether the orbits of the evolution
operator~$T$ generate a frame for a Hilbert space $\Hil$ \cite{second}.
Thus, the basic ingredients are a Hilbert space $\Hil$, a bounded operator $T\colon \Hil \to \Hil$, and a collection of vectors $\{\gamma_j\}_{j=1}^E\subset\Hil$, where $E\in\{1, 2, \ldots, \infty\}$. Typically, the Hilbert space is infinite-dimensional. The fundamental problem is to determine when the family
\begin{equation}\label{uman1}
\mathcal{F}_{T, \gamma}=\bigcup_{j=1}^E\{T^n\gamma_j\}_{n=0}^\infty
\end{equation}
forms a frame for $\Hil$. Such frames are known as \emph{dynamical frames} or, equivalently, \emph{frames of iterations}. The vectors $\gamma_j$ in \eqref{uman1} are called the \emph{generating vectors}, while  the bounded operator ${T}$ is called the \emph{generator operator}. 

An effective approach to studying the family $\mathcal{F}_{T, \gamma}$ is to replace it with an equivalent family $\mathcal{F}_{T', \gamma'}$
that has a simpler structure. The first analysis is performed for $\mathcal{F}_{T', \gamma'}$ and then the obtained results are transferred back to the original family $\mathcal{F}_{T, \gamma}$. We recall  that two families $\mathcal{F}_{T, \gamma}$ and $\mathcal{F}_{T', \gamma'}$ with generators $\{\gamma_j\}_{j=1}^E\subset\Hil$
 and $\{\gamma_j'\}_{j=1}^E\subset\Hil'$, respectively, are called \emph{equivalent} 
if there exists an isomorphism $V \colon \Hil\to\Hil'$ such that
\begin{equation*}
 T'=VTV^{-1} \quad \mbox{and} \quad\gamma_j'=V\gamma_j, \quad j=1, \ldots, E.
\end{equation*}

One natural application of the above method is the following. Let $\mathcal{F}_{T, \gamma}$ be a frame in $\Hil$ with the frame operator $S$. Then the reconstruction formula can be rewritten as
\begin{equation*}
f=\sum_{n=0}^\infty\sum_{j=1}^E\langle{f, S^{-1/2}T^n\gamma_j}\rangle{S^{-1/2}T^n\gamma_j}
=\sum_{n=0}^\infty\sum_{j=1}^E\langle{f, T'^n\gamma_j'}\rangle{T'^n\gamma_j'},
\end{equation*}
where $\gamma_j'=S^{-1/2}\gamma_j$.
The relation obtained means that $\mathcal{F}_{T, \gamma}$ is equivalent to the Parseval frame $\mathcal{F}_{T', \gamma'}$ with $T'=S^{-1/2}TS^{1/2}$. This known result \cite{Rev}  allows us to reduce the investigation of $\mathcal{F}_{T,\gamma}$ to the case of a Parseval dynamical frame generated by a contraction operator. 

Traditionally, an equivalent dynamical frame $\mathcal{F}_{T', \gamma'}$ is constructed in a~suitable functional space. More precisely, the Hilbert space $\Hil'$
is identified with a closed subspace $K$ of the Hardy space $H^2$ of vector-valued analytic functions on the open unit disk in the complex plane.
The subspace $K$ is assumed to be invariant under the adjoint of the unilateral shift on $H^2$, and the operator $T'$ is then identified with the compression of the shift to $K$. This functional model allows us to apply the well-developed tools of harmonic analysis in Hardy spaces \cite{Nik}, leading to many interesting results in dynamical frames \cite{Nie, chri2, 50, 48}. However, in our view, the use of such advanced methods is not always essential since most of the existing results can be obtained by direct arguments without resorting to the machinery of harmonic analysis. Motivated by this observation, we develop a new approach to dynamical frames based on the dilation theory of contractions due to B. Sz.-Nagy and C.~Foiaş \cite{SN}. This framework yields direct operator-theoretic proofs of several known results and, to the best of our knowledge, is the first to highlight the fundamental role of the characteristic function $\Theta_T(\lambda)$ and the isometric dilation of a contraction $T$ in the analysis of dynamical frame generators. Our approach is partly inspired by \cite{Nie},
where minimal isometric dilations were used in the study of dynamical frames. 

The paper is organized as follows. In Section \ref{SecII},
key definitions and results of the Sz.-Nagy–Foiaș dilation
theory \cite{SN} are recalled. The material is presented in 
a way that emphasizes our main objective: to show that
dynamical frames can be studied naturally in the framework of
dilation theory. To highlight this connection, we adopt the
terminology and notation of \cite{SNF, SN} with one
exception. We retain the term strongly stable operator, 
which has become standard in the theory of dynamical frames;
in the Sz.-Nagy–Foiaș terminology, such operators belong to
the class $C_{*0}$.

The purpose of Section \ref{SecIII} is twofold. 
First, we establish several new results on dynamical frames.
Second, we show that a number of recent known results on
dynamical frames follow directly from well-known results in 
dilation theory. 
We begin by recalling the standard definitions 
of frame theory, following the monograph \cite{chri}. 
After briefly examining the structure of the minimal
isometric dilation, in Proposition \ref{aaa1} 
we establish an explicit formula for a~special class
of dynamical frames generated by a contraction. 
This formula was originally obtained 
in \cite[Proposition 3.6]{Nie} using a different method. 
Our proof is closer in spirit to the approach developed in \cite[Lemma~3.7]{Nie}. 
The first result in this direction appeared 
in \cite[Theorem 8]{third}, 
where only the case of Parseval frames was considered. 

A characteristic feature of the dynamical frame \eqref{umc}
considered in Proposition~\ref{aaa1} is that 
all its generating vectors $\{\gamma_j\}_{j=1}^E$ belong to
the defect subspace $\mathfrak{D}_{T^*}$ of contraction~$T$. Motivated by this observation, we say that a dynamical frame $\mathcal{F}_{T,\gamma}$ is \emph{defect-space generated} if all its generating vectors belong to $\mathfrak{D}_{T^*}$. In Theorem \ref{uuu1212}, we show that every Parseval dynamical frame is defect-space generated.

In contrast to the methods used in the existing literature (see, for instance, \cite{Nie, chri2, 48}),we characterize the  generator and the generating vectors 
of a general dynamical frame in terms of the Stein equation (Propositions \ref{uman88}, \ref{uman333}). 
We then examine the frame index $\gamma(T)$, introduced in \cite{Nie}, and determine 
its value in several special cases. The case of jointly strongly stable generators is characterized in terms of the characteristic function $\Theta_T(\lambda)$ (Lemma \ref{uman500}) and the shift index $\chi(T)$ (Corollary~\ref{uman501}). 

Finally, we illustrate the method developed in this paper through simple examples involving contractions $T$ acting diagonally with respect to an orthonormal basis of $\Hil$, as well as weighted unilateral shifts. Our analysis naturally leads to the concept of Carleson frames and highlights the role of the characteristic function (Theorem \ref{the1} and Corollary \ref{the2}).

\section{Selected results of  the Sz.-Nagy–Foiaș dilation theory}\label{SecII}

\subsection{Notation and terminology}\label{sec.2.1}

Throughout the paper, all operators \textit{are assumed to be linear and bounded}. We use the notation $\mathbb{N}=\{1, 2, \ldots\}$, $\mathbb{N}_0=\{0, 1, 2, \ldots\}$,  $\mathbb{D}=\{\lambda\in\mathbb{C} : |\lambda|<1\}$, and $\mathbb{T}=\{\lambda\in\mathbb{C} : |\lambda|=1\}$. 

Let $\Hil$ be a separable infinite-dimensional Hilbert space with inner product $\langle\cdot, \cdot\rangle$, and let 
$T$ be an operator in $\Hil$. We denote its range by $\mathcal{R}(T)$ and its kernel by $\ker T$. The notation $\overline{\mathcal{R}(T)}$ stands for the closure of $\mathcal{R}(T)$ in $\Hil$. The adjoint of $T$ is denoted by~$T^\ast$. 
 An operator $T$ is called \emph{nonnegative} if $\langle Tf,f\rangle \geq 0$ for all $f\in\Hil$, and \emph{uniformly positive} if 
$\langle Tf,f\rangle \geq \alpha\|f\|^2$
for some constant $\alpha>0$ and all $f\in\Hil$. An operator $T$ is called a \emph{contraction} if $\|T\|\leq{1}$. We say that $T$ is \emph{strongly stable} if $\lim_{n\to\infty}T^nf=0$ for every $f\in\Hil$, and \emph{jointly strongly stable} if both $T$ and its adjoint $T^*$ are strongly stable. 
We denote the class of jointly strongly stable contractions by $C_{00}$.

A contraction $T$ is called \emph{completely nonunitary} if it has no nontrivial reducing subspaces on which its restriction is unitary. Every strongly stable contraction is completely nonunitary. A completely nonunitary contraction~$T$ belongs to \emph{the class} $C_0$ if and only if there exists an inner function \mbox{$u\in{H^\infty(\mathbb{D})}$}, where $H^\infty(\mathbb{D})$ is the space of all bounded analytic functions on the unit disk, such that $u(T)=0$.
The class $C_0$  is contained in $C_{00}$. 

Let $\mathcal{L}_1$ and $\mathcal{L}_2$ be closed subspaces of Hilbert spaces  
 ${\mathcal H}_1$ and ${\mathcal H}_2$, respectively. An operator  $V \colon \mathcal{H}_1 \to \mathcal{H}_2$ 
 that maps $\mathcal{L}_1$ isometrically onto  $\mathcal{L}_2$
 and annihilates $\mathcal{H}_1\ominus\mathcal{L}_1$ is called a \emph{partial isometry}. Then $\mathcal{L}_1$ is called the \emph{initial space}, and $\mathcal{L}_2$ is called the \emph{final space} of the partial isometry $V$. 

\subsection{Dilation and characteristic function}

For a contraction $T$ in $\Hil$, the defect operators are defined by
\[
D_T = (I - T^*T)^{1/2}, \quad  
D_{T^*} = (I - TT^*)^{1/2},
\]
and the corresponding defect spaces are
$\mathfrak{D}_T = \overline{\mathcal{R}(D_T)}$ and  
$\mathfrak{D}_{T^*} = \overline{\mathcal{R}(D_{T^*})}$.
The defect indices of $T$ are $\mathfrak{d}_T = \dim \mathfrak{D}_T 
\ \mbox{and} \
\mathfrak{d}_{T^*} = \dim \mathfrak{D}_{T^*}$.
The defect operators satisfy the  intertwining relations:
\begin{equation}\label{uman11}
T D_T = D_{T^*} T, 
\quad 
D_T T^* = T^* D_{T^*}.
\end{equation}
An \emph{isometric dilation} of a contraction $T$ is an isometry $\mathbf{U}$ on a Hilbert space $\mathbf{H}\supset\Hil$, such that
\begin{equation*}
T^n = P_{\Hil} \mathbf{U}^n \big|_{\Hil}, \quad n\in\mathbb{N},
\end{equation*}
where $P_{\Hil}$ is the orthogonal projection of $\mathbf{H}$ onto $\Hil$. The dilation is called \emph{minimal} if
\[
\mathbf{H} = \overline{\operatorname{span}}\{ \mathbf{U}^n h : h \in \Hil,\ n \in\mathbb{N}_0 \}.
\]

The following explicit construction of the minimal isometric dilation of~$T$ is due to \cite{SN}:
\begin{equation}\label{uman7}
\mathbf{U}\mathbf{h} = (T h_0, D_T h_0, h_1, \ldots ), \quad \mathbf{h}\in \mathbf{H},
\end{equation}
where the isometry $\mathbf{U}$ acts on the Hilbert space 
\begin{equation}\label{uman8}
\textbf{H}=\Hil\oplus\bigoplus_{n=1}^{\infty}\mathfrak{D}_T,
\end{equation}
whose vectors are of the form 
\begin{equation*}
\textbf{h}=(h_0,h_1,\ldots ) \textrm{ with } h_0\in\Hil,  h_n\in\mathfrak{D}_T \text{ for } n \in\mathbb{N}, \quad 
\|\textbf{h}\|^2=\sum_{n=0}^{\infty}\|h_n\|^2<\infty.
\end{equation*}

A closed subspace $\mathfrak{L}$ of $\mathbf{H}$ is called a \emph{wandering subspace} for an isometry~$\mathbf{U}$ if $\mathfrak{L}\perp\mathbf{U}^n\mathfrak{L}$ for every  $n\in\mathbb{N}$. An isometry $\mathbf{U}$ on $\mathbf{H}$ is called a \emph{unilateral shift} if there exists a wandering subspace $\mathfrak{L}$ such that
$\bigoplus_{n=0}^\infty\mathbf{U}^n\mathfrak{L}=\mathbf{H}$. The dimension of $\mathfrak{L}$ is called \emph{the multiplicity} of $\bf U$. A unilateral shift with multiplicity $1$ is called \emph{simple}.

The \emph{characteristic function} of a contraction $T$ is the operator-valued function defined on $\mathbb{D}$, with values in the space of bounded linear operators from $\mathfrak{D}_T$ into $\mathfrak{D}_{T^*}$, given by
\begin{equation*}
\begin{split}
\Theta_T(\lambda)&
=[-T+\lambda D_{T^*}(I-\lambda T^*)^{-1}D_T
]\big|_{\mathfrak{D}_T}\\
&=[-T+\lambda\sum_{n=0}^\infty{\lambda^n}D_{T^*}{T^*}^{n}D_T]\big|_{\mathfrak{D}_T},
\quad \lambda\in\mathbb{D}.
\end{split}
\end{equation*}
By construction, the function $\Theta_T$ is analytic on $\mathbb{D}$. The limit $
\Theta_T(e^{it})=\lim_{r\to{1^-}}\Theta_T(re^{it})$
 exists for almost all $e^{it}$ in the sense of $L_2$-convergence \mbox{\cite[Chapter 6, \S 1]{SN}}. 

\subsection{Contractions with strongly stable adjoint}\label{sec.2.3}

Let $T$ be a contraction in $\Hil$ with strongly stable adjoint.
It follows from \cite[Chapter 1, Theorem 3.2]{SN} that $T$ must be completely nonunitary. 
Completely nonunitary contractions with strongly stable adjoints can be characterized in terms of their characteristic functions.

\begin{lem} [{\cite[Chapter 6, \S 3]{SN}}] \label{uman15}
Let $T$ be a completely nonunitary contraction in $\Hil$. The following are equivalent:
\begin{itemize}
\item[(i)] $T$ is a contraction with strongly stable adjoint $T^*$;
\item[(ii)]  the characteristic function $\Theta_T(\lambda)$ is an \emph{inner function}, that is, \linebreak \mbox{$\Theta_T(e^{it}) \colon \mathfrak{D}_T \to \mathfrak{D}_{T^*}$} is an isometric operator a.e. 
\end{itemize}
\end{lem}
The characteristic function $\Theta_T(\lambda)$ determines a contraction $T$ with strongly stable adjoint up to unitary equivalence. The spectrum of $T$ is determined by its characteristic function, and a functional model can be constructed.
\begin{lem}[{\cite[Chapter 6, \S 2, \S 4]{SN}}]\label{uman002}
Let $T$ be a contraction  with strongly stable adjoint. Then:
\begin{itemize}
\item[(i)] $\lambda\in\mathbb{D}$ is an eigenvalue of $T$ if and only if $0$ is an eigenvalue of $\Theta_T(\lambda)$;
\item[(ii)] the operator $T$ is unitarily equivalent to the contraction $\sf T$ acting in the Hilbert space $K=H^2(\mathfrak{D}_{T^*})\ominus\Theta_T{H^2(\mathfrak{D}_{T}})$ and  determined by the formula 
\begin{equation*}
{\sf T}u(\lambda)=P_{K}{\lambda}u(\lambda), \quad {\sf T}^*u(\lambda)=\frac{1}{\lambda}[u(\lambda)-u(0)], \quad u\in{K}.
\end{equation*}
\end{itemize}
\end{lem}

Let $T$ be a contraction with strongly stable adjoint. It follows from Lemma \ref{uman15} that $\mathfrak{d}_{T}\leq\mathfrak{d}_{T^*}$. 
The sign in the last inequality, whether $<$ or $=$, plays an important role in determining the properties of $T$. In particular, if $\mathfrak{d}_{T}<\mathfrak{d}_{T^*}$,
then $T$ has \emph{shift-type invariant subspaces}. More precisely, there exist invariant subspaces $\mathcal{M}_j$ of $T$ such that restrictions $T|_{\mathcal{M}_j}$ are similar to simple unilateral shifts acting on $\mathcal{M}_j$ \cite{Ker}. Another perspective on this subject is provided by the concept of injection of unilateral shifts into contractions.  Let $S$ be a simple unilateral shift $S \colon u \to \lambda{u}$ in ${H}^2(\mathbb{D})$. Denote by $S_m$ the unilateral shift in 
\begin{equation*}
H_m^2(\mathbb{D})={H^2(\mathbb{D})}\oplus{H^2(\mathbb{D})}\oplus\dots\oplus{H^2(\mathbb{D})}, \quad 1\leq{m}\leq\infty
\end{equation*}
defined as the direct sum of $m$ copies of $S$. We say that $S_m$ is \emph{injected into} $T$ if there exists an injective operator $X \colon H_m^2(\mathbb{D})\to \Hil$ such that $XS_m=TX$.

A \emph{shift index} $\chi(T)$ for a contraction $T$ is defined 
\begin{equation*}
\chi(T)=\sup \{m\in\mathbb{N}_{0}\cup\{\infty\} :  S_m \mbox{ injects into } T\}.
\end{equation*}

The following result is an immediate consequence of \cite{SNF} and \cite[Theorem 2]{KT}.

\begin{prop}\label{uman400} 
Let $T$ be a contraction with strongly stable adjoint. The following holds:
\begin{itemize}
\item[(i)]
if $\mathfrak{d}_{T}<\mathfrak{d}_{T^*}$, then 
$\chi(T)=\mathfrak{d}_{T^*}-\mathfrak{d}_{T}\not=0$;
\item[(ii)] if $T$ is jointly strongly stable (i.e. $T\in{C_{00}}$), then $\chi(T)=0$ when $T\in C_0$, whereas $\chi(T)=\infty$ when $T\notin C_0$.
\end{itemize}
\end{prop}

\section{Dynamical frames}\label{SecIII}

\subsection{Standard definitions}\label{sec.3.1}

A family of vectors $\{\varphi_j\}_{j=0}^\infty\subset\Hil$ is called a \emph{frame} for a Hilbert space $\Hil$ if there exist constants $A$ and $B$ satisfying $0<A\leq B<\infty$ such that
\begin{equation}\label{uman63}
A\|f\|^2\leq\sum_{j=0}^\infty{|\langle f, \varphi_j\rangle |^2}\leq{B}\|f\|^2, \quad \forall{f}\in \Hil.
\end{equation}
A family of vectors $\{\varphi_j\}_{j=0}^\infty$ is called a \emph{frame sequence} if it is a frame for $\overline{\operatorname{span}}\{\varphi_j\}_{j=0}^\infty$.
An \emph{$A$-tight frame} is a frame $\{\varphi_j\}_{j=0}^\infty$
for which the frame bounds in \eqref{uman63} can be chosen to be $A=B$. A $1$-tight frame is called \emph{a~Parseval frame}. 

For a frame $\{\varphi_j\}_{j=0}^\infty$, the corresponding
\emph{frame operator} $S$, defined by $S{f}=\sum_{j=0}^\infty\langle f, \varphi_j\rangle \varphi_j$, is uniformly positive in $\Hil$. 
The frame operators for $A$-tight frames coincide with $AI$, 
where $I$ is the identity operator in $\Hil$. If $S$ is a frame 
operator for a frame  $\{\varphi_j\}_{j=0}^\infty$, then the set 
$\{\tau_j=S^{-1/2}\varphi_j\}_{j=0}^\infty$ is a~Parseval frame.

A sequence of vectors $\{\varphi_j\}_{j=0}^E$, where $E\in\mathbb{N}\cup\{\infty\}$, is called a \emph{Bessel sequence} if there exist a constant $B>0$ such that
\begin{equation*}
\sum_{j=0}^E{|\langle f, \varphi_j\rangle |^2}\leq{B}\|f\|^2, \quad f \in \Hil.
\end{equation*}
Each Bessel sequence $\{\varphi_j\}_{j=0}^E$ determines a bounded
operator $\theta \colon \Hil\to\ell_2(J)$, where $J=\{0, 1, \ldots, E\}$, 
\begin{equation}\label{K31}
\theta f=\{\langle f,\varphi_j\rangle\}_{j=0}^E,
\quad f\in\Hil.
\end{equation}
called the \emph{analysis operator}.
Its adjoint $\theta^\ast \colon \ell_2(J) \to \Hil$, given by
\begin{equation*}
\theta^*\{c_j\}_{j=0}^{E}=\sum_{j=0}^E{c_j}\varphi_j, \quad \{c_j\}_{j=0}^{E}\in\ell_2(J),
\end{equation*}
is a bounded operator and is called the \emph{synthesis operator}. The corresponding reconstruction operator $S$, defined as 
\begin{equation*}
Sf=\theta^*\theta{f}=\sum_{j=0}^E\langle f, \varphi_j\rangle \varphi_j, \quad f \in \Hil,
\end{equation*}
is nonnegative on $\Hil$. 

\subsection{Dynamical frames constructed with the use of unilateral shift}

Let $\mathbf{U}$ be a minimal isometric dilation of a contraction $T$ defined by \eqref{uman7} and \eqref{uman8}. The~next result follows directly from \cite[Theorem II.2.1]{SN}.

\begin{prop}\label{uman9}
The  contraction $T^*$ is strongly stable
if and only if  the isometric dilation $\mathbf{U}$ of $T$ is a unilateral shift in $\mathbf{H}$ with the wandering subspace
\begin{equation}\label{uman10}
\mathfrak{L}_*
= \overline{(I-\mathbf{U}T^*)\Hil} 
\end{equation}
\end{prop}

The subspace $\mathfrak{L}_*$ admits a simple characterization. 

\begin{lem}\label{p1p}
The subspace $\mathfrak{L}_*$ has the form
\begin{equation*}
\mathfrak{L}_*=\{(D_{T^*}\tau, -T^*\tau, 0, \ldots) : \tau\in\mathfrak{D}_{T^*} \}
\end{equation*}
and the operator $W \colon \mathfrak{D}_{T^*}\to\mathfrak{L}_*$ defined by
\begin{equation}\label{uman123}
W\tau=(D_{T^*}\tau, -T^*\tau, 0, \ldots), \quad \tau\in\mathfrak{D}_{T^*},
\end{equation}
is unitary.
\end{lem}

\begin{proof} 
In view of \eqref{uman11}--\eqref{uman8} and \eqref{uman10}, we have
\begin{equation*}
\begin{split}
\mathfrak{L}_*=\overline{(I-\mathbf{U}T^*)\Hil}
&=\overline{\{(I-TT^*)h_0,
-D_TT^*h_0, 0, \ldots) : h_0\in\Hil\}}\\
&=\overline{\{(D_{T^*}^2h_0,
-T^*D_{T^*}h_0, 0, \ldots) : h_0\in\Hil\}}\\ 
&=\overline{\{(D_{T^*}\tau,
-T^*\tau, 0, \ldots) : \tau\in\mathcal{R}(D_{T^*})\}}\\ 
&=\overline{\{(D_{T^*}\tau,
-T^*\tau, 0, \ldots) : \tau\in\mathfrak{D}_{T^*}\}}=\overline{\mathcal{R}(W)},
\end{split}
\end{equation*}
where $W$ is defined by \eqref{uman123}. Denote
$\mathbf{h}=(D_{T^*}{\tau}, -T^*{\tau}, 0, \ldots)$, where $\tau\in\mathfrak{D}_{T^*}$.  
It follows from
\begin{equation*}
\langle{W\tau}, W\tau\rangle=\langle\mathbf{h}, \mathbf{h}\rangle=\langle(I-TT^*)\tau, \tau\rangle+\langle TT^*\tau, \tau\rangle=\langle\tau, \tau\rangle, \quad \tau\in\mathfrak{D}_{T^*}
\end{equation*}
that the operator $W$ is a unitary mapping of $\mathfrak{D}_{T^*}$ onto
\begin{equation*}
\mathfrak{L}_*=\overline{\mathcal{R}(W)}=\mathcal{R}(W)=\{(D_{T^*}\tau, -T^*\tau, 0, \ldots) : \tau\in\mathfrak{D}_{T^*} \}.
\end{equation*}
This completes the proof.
\end{proof}

Consider a frame $\{\tau_j\}_{j=1}^E$ $(E\in\mathbb{N}\cup\{\infty\})$ for the Hilbert space $\mathfrak{D}_{T^*}$ with frame bounds $A\leq{B}$, and
define
\begin{equation*}
W\tau_j=\mathbf{e}_j=(D_{T^*}{\tau_j}, -T^*{\tau_j}, 0, \ldots), \quad j\in\{1, 2, \ldots, E\}.  
\end{equation*} 
By Lemma \ref{p1p}, $W$ is unitary. Hence, the collection of vectors $\{\mathbf{e}_j\}_{j=1}^E$ is a~frame for the space $\mathfrak{L}_*$ with the same frame bounds $A$ and $B$.  
In view of~\eqref{uman7},
the vectors 
\begin{equation*}
\mathbf{e}_j^1=\mathbf{U}\mathbf{e}_j=(TD_{T^*}{\tau_j},D_TD_{T^*}\tau_j, -T^*{\tau_j}, 0, \ldots),   
\end{equation*}
where $j=1, 2, \ldots, E$, form a frame for $\mathbf{U}\mathfrak{L}_*$. Similarly, for every $n\geq2$, the~collection of vectors
\begin{equation*}
\mathbf{e}_j^n=\mathbf{U}^n\mathbf{e}_{j}
=(T^nD_{T^*}{\tau_j}, D_TT^{n-1}D_{T^*}{\tau_j},\ldots, D_TD_{T^*}\tau_j, -T^*{\tau_j}, 0, \ldots), 
\end{equation*}
where $j=1, 2, \ldots, E$, gives rise to a frame for
$\mathbf{U}^n\mathfrak{L}_*$.

\begin{prop}\label{aaa1}
Assume that $T$ is a contraction with a strongly stable adjoint.  Then, for an arbitrary frame $\{\tau_j\}_{j=1}^{E}$ for $\mathfrak{D}_{T^*}$, the set 
\begin{equation}\label{umc}
\mathcal{F}_{T,\gamma}=\bigcup_{j=1}^{E}\{T^n\gamma_j\}_{n=0}^\infty, \quad \gamma_j=D_{T^*}\tau_j, \quad j=1, 2, \ldots, E,
\end{equation} 
is a dynamical frame for $\Hil$  with the same frame bounds as the frame $\{\tau_j\}_{j=1}^E$ in $\mathfrak{D}_{T^*}$. The formula \eqref{umc} defines an $A$-tight dynamical frame whenever $\{\tau_j\}_{j=1}^E$ is an $A$-tight frame in~$\mathfrak{D}_{T^*}$.
\end{prop}

\begin{proof} 
It follows from Proposition \ref{uman9} that
\begin{equation*}
\mathbf{H}=\bigoplus_{n=0}^{\infty}\mathbf{U}^n\mathfrak{L}_*
\end{equation*}
Therefore, the vector collection \begin{equation}\label{uman13}
\bigoplus_{n=0}^\infty\{{\bf U}^n{\bf e}_j\}_{j=1}^E=\bigcup_{j=1}^E\{{\bf U}^n{\bf e}_j\}_{n=0}^\infty=\mathcal{F}_{{\bf U}, {\bf e}}
\end{equation}
 is a dynamical frame of $\mathbf{H}$ with the same frame bounds  as the frame $\{\tau_j\}_{j=1}^E$. In particular, \eqref{uman13} turns out to be an $A$-tight  frame whenever $\{\tau_j\}_{j=1}^E$ is an $A$-tight frame in $\mathfrak{D}_{T^*}$.
By~projecting \eqref{uman13} onto $\Hil$, we obtain the ($A$-tight) frame \eqref{umc}, which completes the proof.
\end{proof}

As follows from Proposition \ref{aaa1}, a dynamical frame \eqref{umc} with contraction generator $T$ and generating vectors $\{\gamma_j\}_{j=1}^E$ has the characteristic property that all $\gamma_j\in\mathfrak{D}_{T^*}$. For this reason, we say that a dynamical frame $\mathcal{F}_{T, \gamma}$ is \emph{defect-space generated} if all its generating vectors belong to the defect subspace $\mathfrak{D}_{T^*}$ of $T$.

\subsection{Parseval dynamical frames}\label{secIII.2}

Let $\mathcal{F}_{{T}, \gamma_j}$ be an arbitrary Parseval dynamical  frame defined by \eqref{uman1}. Then
\begin{equation}\label{uuu2b}
(I-T^k{T^*}^k)f=\sum_{n=0}^{k-1}\sum_{j=1}^E\langle{f, T^n\gamma_j}\rangle{T^n\gamma_j}
\end{equation}
and
\begin{equation}\label{uman6}
\|f\|^2-\|{T^*}^kf\|^2=\sum_{n=0}^{k-1}\sum_{j=1}^E|\langle{f, T^n\gamma_j}\rangle|^2
\end{equation}

Taking $k=1$ in \eqref{uman6} gives that $T^*$ is a contraction on $\Hil$. Since \mbox{$\|T\|=\|T^*\|$}, the~operator $T$ is also a contraction. Furthermore, for all $f\in\Hil$,
\begin{equation*}
\lim_{k\to\infty}\|{T^*}^kf\|^2=\|f\|^2-\lim_{k\to\infty}\sum_{n=0}^{k-1}\sum_{j=1}^E|\langle{f, T^n\gamma_j}\rangle|^2=\|f\|^2-\|f\|^2=0,
\end{equation*}
implying that $T^*$ is strongly stable. Thus, every Parseval dynamical frame is generated by a~contraction whose adjoint is strongly stable. This result was first proved in \cite{third} by a different method. The converse statement is also true.
 
\begin{prop}\label{uman14}
Let $T$ be a contraction in $\Hil$. The following are equivalent:
\begin{itemize}
\item[(i)] there is a collection of generators $\{\gamma_j\}$ such that $\mathcal{F}_{T, \gamma}$ is a Parseval dynamical frame;
\item[(ii)] the contraction $T^*$ is strongly stable;
\item[(iii)] a minimal isometric dilation $\mathbf{U}$ of $T$ is a unilateral shift in $\mathbf{H}$;
\item[(iv)] $T$ is completely nonunitary  and its  characteristic function $\Theta_T(\lambda)$ is an inner function.
\end{itemize}
\end{prop}

\begin{proof} 
The implication (i)$\Rightarrow$(ii) has been established immediately before Pro\-po\-sition \ref{uman14}. The implication (ii)$\Rightarrow$(i) follows from Proposition \ref{aaa1}. The equivalence of (ii) and (iii) follows from Proposition~\ref{uman9}, while the equivalence of (ii) and (iv) follows from Lemma~\ref{uman15}.
\end{proof}

\begin{rem}
The implication (ii)$\Rightarrow$(i)  was first proved in \cite[Theorem 8]{third}. The~equivalence between (i) and (ii) was established in \cite{Nie} using the functional models (equivalent frames) mentioned in the Introduction.
\end{rem} 

\begin{thm}\label{uuu1212}
Each Parseval dynamical frame $\mathcal{F}_{T,\gamma}$ is defect space generated and $\mathcal{F}_{T,\gamma}$ is defined by the formula \eqref{umc} for a suitable choice of a Parseval frame $\{\tau_j\}_{j=1}^E$ in $\mathfrak{D}_{T^*}$.
\end{thm}

\begin{proof} 
As follows from \eqref{uuu2b} with $k=1$,
\begin{equation}\label{uuu4b}
D_{T^*}^2f=(I-T{T^*})f=\sum_{j=1}^E\langle{f, \gamma_j}\rangle{\gamma_j}.
\end{equation}
Since $\mathcal{F}_{T,\gamma}$ is a frame for $\Hil$ defined by \eqref{uman1}, the sequence
$\{\gamma_j\}_{j=1}^E$ is a Bessel sequence. In view of \eqref{uuu4b}, $D^2_{T^*}$~coincides with the reconstruction operator of $\{\gamma_j\}_{j=1}^E$ and has the form $D_{T^*}^2=\theta^*\theta$, where $\theta$ is the analysis operator of $\{\gamma_j\}_{j=1}^E$ (see \eqref{K31}). This means that the polar decomposition of $\theta$ is \mbox{$\theta=VD_{T^*}$}, where $V \colon \Hil \to \ell_2(J)$, $J=\{1,2,\ldots, E\}$, is a partial isometry with initial space $\mathfrak{D}_{T^*}$ and final space $\overline{\mathcal{R}(\theta)}$. 

Denote $\tau_j=V^*\delta_j$, where $\{\delta_j\}_{j=1}^E$ is the canonical basis of $\ell_2(J)$. By construction, $\{\tau_j\}_{j=1}^E\subset\mathfrak{D}_{T^*}$ and
$D_{T^*}\tau_j=D_{T^*}V^*\delta_j=\theta^*\delta_j=\gamma_j$. Furthermore, $\{\tau_j\}_{j=1}^E$ is a Parseval frame in $\mathfrak{D}_{T^*}$ because 
\begin{equation*}
\sum_{j=1}^E|\langle f, \tau_j\rangle|^2=\sum_{j=1}^E|\langle Vf, \delta_j\rangle|^2=\|Vf\|^2_{\ell_2(J)}=\|f\|^2, \quad f\in\mathfrak{D}_{T^*}. 
\end{equation*}
From Proposition \ref{aaa1} it follows that $\mathcal{F}_{T, \gamma}$ is defined by the formula \eqref{umc}.
\end{proof}

The following result characterizes the generating vectors of Parseval dynamical frames and complements Proposition \ref{uman14}, which describes the properties of the generator operators.

\begin{cor}\label{ccc21}
Let  $\mathcal{F}_{T,\gamma}$ be a Parseval  dynamical frame. Then
\begin{equation}\label{bbb23}
\mathfrak{D}_{T^*}=\overline{\operatorname{span}}\{\gamma_j\}_{j=1}^E.
\end{equation}
The set of generating vectors $\{\gamma_j\}_{j=1}^E$ is a frame sequence in $\Hil$ if and only if \mbox{$\mathcal{R}(D_{T^*})=\overline{\mathcal{R}(D_{T^*})}$}. If $E$ in \eqref{bbb23} is finite, then 
\begin{equation}\label{uuu166}
\sum_{j=1}^E\|\gamma_j\|^2=\mathfrak{d}_{T^*}-\operatorname{Tr} \ TT^*|_{\mathfrak{D}_{T^*}},
\end{equation}
where $TT^*|_{\mathfrak{D}_{T^*}}$ denotes the restriction of $TT^*$ to the reducing subspace $\mathfrak{D}_{T^*}$. 
\end{cor}

\begin{proof}
Formula \eqref{bbb23} follows immediately from the proof of Theorem \ref{uuu1212}, since 
$\gamma_j=D_{T^*}\tau_j$, formula \eqref{uuu4b} holds, and $\{\tau_j\}_{j=1}^E$ is a Parseval frame in~$\mathfrak{D}_{T^*}$. 

Assume that $\{\gamma_j\}_{j=1}^E$ is a frame sequence.  Then, in view of \eqref{bbb23}, $\{\gamma_j\}_{j=1}^E$ is a frame for $\mathfrak{D}_{T^*}$. By means of \eqref{uuu4b}, the restriction of $D_{T^*}^2$ to~$\mathfrak{D}_{T^*}$ coincides with the frame operator $S$ of $\{\gamma_j\}_{j=1}^E$. 
Therefore\footnote{Note that all operators are considered in $\mathfrak{D}_{T^*}$}, $0\in\rho(D_{T^*}^2)$ and $\mathcal{R}(D_{T^*})=\overline{\mathcal{R}(D_{T^*})}$.

Conversely, the condition 
${\mathcal{R}(D_{T^*})}=\overline{\mathcal{R}(D_{T^*})}$ means that $0\in\rho(D_{T^*})$ and, hence, $\{\gamma_j\}_{j=1}^E$ is a frame sequence in $\Hil$ since $\{\tau_{j}\}_{j=1}^E$ is a Parseval frame in $\mathfrak{D}_{T^*}$.

Taking $\gamma_j=D_{T^*}\tau_j$ into account, we obtain
\[
\langle\gamma_j, \gamma_j\rangle=\langle(I-TT^*)\tau_j, \tau_j\rangle=\|\tau_j\|^2-\langle{TT^*}\tau_j, \tau_j\rangle.
\]
Summing this relation over $j$ yields \eqref{uuu166}, as
$\sum_{j=1}^E \|\tau_j\|^2
=\dim \mathfrak{D}_{T^*}
=\mathfrak{d}_{T^*}$, while, by~\mbox{\cite[Proposition 4.11]{CL}}, $\sum_{j=1}^E\langle{TT^*}\tau_j, \tau_j\rangle=\operatorname{Tr} \ TT^*|_{\mathfrak{D}_{T^*}}$. 
\end{proof}

\subsection{General dynamical frames}

Similarly to Section \ref{secIII.2}, we first analyze the
properties of generator operators.

\begin{prop}\label{uman88}
Let $T$ be a bounded operator in $\Hil$. The following are equivalent:
\begin{itemize}
\item[(i)] $T$ is a generator  of a dynamical frame $\mathcal{F}_{T,\gamma}$ (see  \eqref{uman1}) with the frame operator $S$;
\item[(ii)] the adjoint $T^*$ is strongly stable and there exists a uniformly positive operator $S$ such that the operator $S- TST^*$ is nonnegative in $\Hil$.
\end{itemize}
\end{prop}

\begin{proof}
Define
\begin{equation}\label{uman77}
T'=S^{-1/2}{T}S^{1/2}.
\end{equation}
Then 
\begin{equation}\label{uman101}
I-T'T'^*=S^{-1/2}(S-TST^*)S^{-1/2}.
\end{equation}
This means that $T'$ is a contraction if and only if $S-TST^*$ is a nonnegative operator. If the condition of (ii) holds, then $T'$ is a contraction with strongly stable adjoint. By Proposition \ref{uman14}  and Theorem \ref{uuu1212}, there exists a Parseval dynamical frame $\mathcal{F}_{T', \gamma'}$, where $\{\gamma_j'=D_{T'^*}\tau_j'\}_{j=1}^E$ and $\{\tau_j'\}_{j=1}^E$ is a Parseval frame in the Hilbert space $\mathfrak{D}_{T'^*}$. By an elementary calculation, using \eqref{uman77}, we deduce that
\begin{equation}\label{uman99}
\mathcal{F}_{T, \gamma}=S^{1/2}\mathcal{F}_{T', \gamma'}, \quad \gamma_j=S^{1/2}\gamma_j', \quad j=1,\ldots,E,
\end{equation}
is a dynamical frame with the generator operator $T$ and the frame operator~$S$. The implication (ii)$\Rightarrow$(i) is proved.

To prove the converse implication (i)$\Rightarrow$(ii), it suffices to observe that 
\begin{equation*}
(S- T^kS{T^*}^k)f=\sum_{n=0}^{k-1}\sum_{j=1}^E \langle f,T^n\gamma_j\rangle{T}^n\gamma_j.
\end{equation*}
This formula guarantees the nonnegativity of $S-T^kS{T^*}^k$ for all $k\in\mathbb{N}$ and yields 
\begin{equation*}
\langle{S^{1/2}{T^*}^kf, S^{1/2}{T^*}^kf}\rangle=\langle{Sf,f}\rangle-\sum_{n=0}^{k-1}\sum_{j=1}^E |\langle f,T^n\gamma_j\rangle|^2.
\end{equation*}
Letting $k\to\infty$ in the last relation, we conclude
$\lim_{k\to\infty}\|S^{1/2}{T^*}^kf\|^2=0$. Hence, $T^*$ is strongly stable.
\end{proof}

\begin{rem}
Another characteristic property of dynamical frame generator operators was established in \cite{Nie, 48}.
\end{rem}

An operator $T$ acting on $\Hil$ is called \emph{$\rho$-contraction} $(\rho>0)$ if there is a~unitary operator 
$\mathbf{U}$ on a Hilbert space $\mathbf{H}\supset\Hil$, such that
$T^n = {\rho}P_{\Hil} \mathbf{U}^n \big|_{\Hil}$,  $n\in\mathbb{N}$, where $P_{\Hil}$ is the orthogonal projection of $\mathbf{H}$ onto $\Hil$. When $\rho=1$ or $\rho=2$, a $\rho$-contraction $T$ is, respectively, a contraction or a numerical contraction, i.e., an operator satisfying
$\sup\{|\langle Tx, x\rangle| : x \in\Hil, \|x\| = 1\} \leq 1
$; see \cite{SN, WG}. 

\begin{cor}
Each $\rho$-contraction $T$ with strongly stable adjoint is a generator operator for a dynamical frame $\mathcal{F}_{T, \gamma}$.
\end{cor}

\begin{proof}
According to \cite[Ch. II, Theorem 8.1]{SN}, every $\rho$-contraction $T$ is similar to a contraction; that is, there exists a boundedly invertible operator $X$ such that $C=X^{-1}TX$
is a contraction on $\Hil$. Denote  $S=XX^*$. The operator $S$ is uniformly positive by construction. Moreover,
\begin{align*}
S-TST^*=X(I-CC^*)X^*\geq{0}.
\end{align*}
The proof is completed by applying Proposition \ref{uman88}.
\end{proof}

\begin{cor}
Let $\mathcal{F}_{T, \gamma}$ be an
$A$-tight dynamical frame. Then $T$ is a contraction and $\mathcal{F}_{T, \gamma}$ is the defect space generated. 
\end{cor}

\begin{proof}
The class of $A$-tight frames is uniquely determined by its frame operator $S=AI$.
Let $\mathcal{F}_{T, \gamma}$ be an
$A$-tight dynamical frame.
 It follows from Proposition \ref{uman88}  that 
\begin{equation*}
S-TST^*=A(I-TT^*) \geq 0.
\end{equation*}
Therefore, $T$ is a contraction. Using \eqref{uman99}, we get
$\mathcal{F}_{T, \gamma}=A^{1/2}\mathcal{F}_{T, \gamma'}$, where $\gamma_j'= A^{-1/2}\gamma_j$ and $\mathcal{F}_{T, \gamma'}$ is a Parseval dynamical frame. The proof is completed by applying Theorem \ref{uuu1212}. 
\end{proof}

We conclude this section with a weakened version of Proposition  \ref{uman88}. 

\begin{prop}\label{uman333}
Assume that a nonnegative  operator $S$ is a solution of the Stein equation\footnote{See \cite{BBF} for details on  the Stein equation}
\begin{equation*}
S-TST^*=P_\gamma, \quad P_\gamma=\sum_{j=1}^E\langle\cdot, \gamma_j\rangle\gamma_j, \quad E\in\mathbb{N},
\end{equation*}
where  $\{\gamma_j\}_{1}^E$ is a family of vectors in $\Hil$  such that the operator $P_\gamma$ is bounded, and $T$ is a bounded operator. Then the set $\mathcal{F}=\bigcup_{j=1}^E\{T^n\gamma_j\}_{n=0}^\infty$ is a Bessel sequence in $\Hil$. Moreover, the set $\mathcal{F}$ is a dynamical frame $\mathcal{F}_{T,\gamma}$ provided that $S$ is uniformly positive and $\lim_{n\to\infty}S^{1/2}{T^*}^nf=0$ for all $f\in\Hil$. In this case, $S$ is the frame operator for $\mathcal{F}_{T,\gamma}$ and $\{\gamma_j\}_{j=1}^E$ are generating vectors.
\end{prop}

\begin{proof} 
Let $S$ be a solution of the Stein equation. Iterating gives
\begin{equation}\label{uman345}
S=\sum_{k=0}^{n-1}T^kP_\gamma{{T^*}^k}+ T^nS{T^*}^n, \quad n\in\mathbb{N}.
\end{equation}
Since $S$ is nonnegative, 
\begin{equation*}
\begin{split}
\langle{Sf, f\rangle} & \geq\sum_{k=0}^{n-1}\langle{T}^kP_\gamma{{T^*}^k}f, f\rangle \\
&=\sum_{k=0}^{n-1}\sum_{j=1}^E|\langle{{T^*}^k}f, \gamma_j\rangle|^2=\sum_{k=0}^{n-1}\sum_{j=1}^E|\langle{f}, T^k\gamma_j\rangle|^2, \quad f\in\Hil. 
\end{split}
\end{equation*}
Hence, $\mathcal{F}=\bigcup_{j=1}^E\{T^n\gamma_j\}_{n=0}^\infty$ is a Bessel sequence in $\Hil$.  Taking the limit as $n\to\infty$ in \eqref{uman345}, we obtain
\begin{equation*}
\begin{split}
\langle{Sf,f}\rangle&=\langle{S'f, f\rangle}+\lim_{n\to\infty}\langle{T}^nS{T^*}^nf, f\rangle\\
&=\sum_{k=0}^{\infty}\sum_{j=1}^E|\langle{f}, T^k\gamma_j\rangle|^2+\lim_{n\to\infty}\|S^{1/2}{T^*}^nf\|^2,
\end{split}
\end{equation*}
where $S'f=\sum_{k=0}^{\infty}\sum_{j=1}^E\langle{f, T^k\gamma_j\rangle}{T^k\gamma_j}$ is the reconstruction operator for the Bessel sequence $\mathcal{F}$. 
If $S$ is uniformly positive and $\lim_{n\to\infty}S^{1/2}{T^*}^nf=0$, then there exists $c>0$ such that 
\begin{align*}
\langle Sf,f\rangle=\sum_{k=0}^{\infty}\sum_{j=1}^E|\langle{f}, T^k\gamma_j\rangle|^2\geq{c\|f\|^2}.
\end{align*}
The last relation means that the Bessel sequence $\mathcal{F}$ is  a frame $\mathcal{F}_{T,\gamma}$ and $S$ is the frame operator for $\mathcal{F}_{T,\gamma}$. 
\end{proof}

\subsection{Frame index $\gamma(T)$}

Let an operator $T$ be a generator of a dynamical frame $\mathcal{F}_{T, \gamma}$. In general, this frame is not uniquely determined by~$T$. There may exist \emph{different sets of generating vectors $\{\gamma_j\}_{j=1}^E$}, possibly with different~$E$, giving rise to distinct frames with the same generator operator $T$.

Following \cite[Definition 4.1]{Nie}, we call the minimum number $E$ of generating vectors $\{\gamma_j\}_{j=1}^E$ of dynamical frames $\mathcal{F}_{T, \gamma}$ with the generator operator~$T$, the \emph{frame index} of~$T$, denoted by $\gamma(T)$.

To find $\gamma(T)$, consider the frame operator $S$ of one of such frames $\mathcal{F}_{T, \gamma}$. It~follows from \eqref{uman99}  and  Theorem \ref{uuu1212} that the minimum value of $E$ is attained when the Parseval frame $\{\tau_j\}_{j=1}^E$ in the proof of Theorem \ref{uuu1212} is~chosen as an orthonormal basis of $\mathfrak{D}_{T'^*}$. Therefore, for a given $T$ and $S$ the minimum of generating vectors $E$ coincides with the defect index $\mathfrak{d}_{T'^*}$ of the contraction~$T'$ defined by \eqref{uman77}. Moreover, in view of \eqref{uman101}, 
\begin{equation*}
\mathfrak{d}_{T'^*}=\dim\overline{\mathcal{R}(S-TST^*)}.
\end{equation*}
Note that this formula is obtained for a given   dynamical frame $\mathcal{F}_{T, \gamma}$ with the frame operator $S$ and therefore $\gamma(T)\leq\mathfrak{d}_{T'^*}$. In view of Proposition \ref{uman88},  the class of all dynamical frames with a given generator $T$ is completely characterized by the set of uniformly positive operators $S$ satisfying \mbox{$S-TST^{*}\geq{0}$}. Consequently, the minimum number of vectors $\{\gamma_j\}_{j=1}^E$ needed to produce a~dynamical frame for a given operator~$T$ can be calculated by the formula
\begin{equation}\label{uman100}
\gamma(T)=\operatorname{min}\dim\overline{\mathcal{R}(S-TST^*)},
\end{equation}
where the minimum is taken over all uniformly positive operators $S$ in $\Hil$ such that $S-TST^*$ is nonnegative.

\begin{rem}
A formula similar to \eqref{uman100} was established in \cite[Theorem~1.3]{Nie} using different methods based on the investigation of equivalent frames defined in a model space of vector-valued analytical functions.
\end{rem}

Let $T$ be a contraction in $\Hil$. Since the identity operator $I$ is an admissible choice of $S$ in \eqref{uman100}, it follows that the frame index  $\gamma(T)$ is bounded above by the defect index~$\mathfrak{d}_{T^*}$, that is,
\begin{equation}\label{uman131}
\gamma(T)\leq\mathfrak{d}_{T^*}.
\end{equation}
In general, the bound in \eqref{uman131} does not need to be reached. We provide several sufficient conditions below under which the equality $\gamma(T)=\mathfrak{d}_{T^*}$ holds.

The decomposition $\Hil=\ker{T^*}\oplus\overline{\mathcal{R}(T)}$ shows that 
\begin{equation*}
\begin{split}
\mathfrak{d}_{T^*}&=\dim\overline{\mathcal{R}(D_{T^*})}=\dim\overline{\mathcal{R}(I-TT^*)}\\
&=\dim\ker{T^*}+\dim\overline{\mathcal{R}(T(I-T^*T))}\geq\dim\ker{T^*}.
\end{split}
\end{equation*}
We first consider the case where the minimum possible value $\mathfrak{d}_{T^*}$ is attained, i.e., the case where $\mathfrak{d}_{T^*}=\dim\ker{T^*}$.

\begin{prop}\label{uman140}
Let $T$ be a contraction with strongly stable adjoint and such that \mbox{$\mathfrak{d}_{T^*}=\dim\ker{T^*}$}. Then $\gamma(T)=\mathfrak{d}_{T^*}$.
\end{prop}
\begin{proof}
 Let
$\mathcal{F}_{T,\gamma}$ be an arbitrary dynamical frame with the generator operator $T$ in $\Hil$. Denote by $S$ and $P_{\mathcal M}$, respectively, the frame operator of $\mathcal{F}_{T, \gamma}$ and the orthogonal projection of $\Hil$ onto  $\mathcal M=\ker T^*$. Since $T^*|_{\mathcal M}=0$, we have
$P_{\mathcal M}(S-TST^*)P_{\mathcal M}
=P_{\mathcal M}SP_{\mathcal M}$.
As $S\geq cI$ for some $c>0$,
\begin{align*}
P_{\mathcal M}SP_{\mathcal M}\geq cP_{\mathcal M}.
\end{align*}
Hence, $\dim\overline{\mathcal{R}(P_{\mathcal M}(S-TST^*)P_{\mathcal M})}=\dim\mathcal M=\dim\ker T^*$. Consequently,
\begin{align*}
\dim\overline{\mathcal{R}(S-TST^*)}\geq \dim\ker T^*=\mathfrak{d}_{T^*}.
\end{align*}
Using \eqref{uman100}, we conclude that 
$\gamma(T)\geq \mathfrak{d}_{T^*}$. Therefore, by~\eqref{uman131}, we obtain $\gamma(T)=\mathfrak{d}_{T^*}$.
\end{proof}

\begin{prop}\label{uman141}
If $\|T\|<1$, then $\gamma(T)=\mathfrak{d}_{T^*}=\infty$.
 \end{prop}
\begin{proof}
If $\|T\|<1$, then $T$ is a contraction with strongly stable adjoint.
Suppose that there exists a dynamical frame $\mathcal{F}_{T,\gamma}$  with a finite number $E$ of
generating vectors $\{\gamma_j\}_{j=1}^E$. Since $\|T\|<1$, we have
\[
\sum_{n=0}^{\infty}\sum_{j=1}^E
\|T^n\gamma_j\|^2
\leq
\frac{1}{1-\|T\|^2}
\sum_{j=1}^E\|\gamma_j\|^2<\infty.
\]
Hence, the frame operator of $\mathcal{F}_{T, \gamma}$, given by
\[
Sf
=
\sum_{n=0}^{\infty}\sum_{j=1}^E
\langle f,T^n\gamma_j\rangle T^n\gamma_j,
\quad f\in\Hil,
\]
is trace class and, in particular, compact. This is impossible for the
frame operator in an infinite-dimensional Hilbert space,
since such an operator is boundedly invertible. Therefore, $E$ must be infinite, and hence, taking into account the definition of $\gamma(T)$ and \eqref{uman131}, we obtain
$\gamma(T)=\mathfrak{d}_{T^*}=\infty$.
\end{proof}

The following corollary is a direct consequence of the proof of Proposition~\ref{uman141}.

\begin{cor}\label{cor20}
If an operator $T$ has a finite frame index $\gamma(T)<\infty$,  then $\|T\|\geq{1}$.
\end{cor}

\begin{rem}
A statement similar to that of Corollary~\ref{cor20} was established in \cite[Theorem 9]{third} by different methods.
\end{rem}

\subsection{Jointly strongly stable generators}
In \cite{Nie}, the jointly strongly stable generators were studied in the setting of vector-valued Hardy spaces. Exploiting the structure of functional models, the authors established a connection between joint  strong stability and the concept of full-range shift-invariant subspaces of Hardy spaces developed in~\cite{Helson}. In contrast to this approach, we characterize the class of contractions~$C_{00}$ using  methods of the Sz.-Nagy-Foiaş dilation theory \cite{SN}. The following result is a direct consequence of Lemma \ref{uman15}.

\begin{lem}\label{uman500}
Let $T$ be a completely nonunitary contraction. Then  $T$ belongs to $C_{00}$ if and only if the characteristic function $\Theta_T(\lambda)$  is inner from both sides, that is, if $\Theta_T(e^{it}) \colon \mathfrak{D}_T \to \mathfrak{D}_{T^*}$ is unitary a.e.
\end{lem}

In the case of finite defect indices, the shift index of contraction (see Section \ref{sec.2.3}) completely characterizes the class $C_{00}$.
\begin{cor}\label{uman501}
Let $T$ be a contraction with strongly stable adjoint and such that $\mathfrak{d}_{T^*}<\infty$. Then $T\in{C_{00}}$ if and only if $\chi(T)=0$.
\end{cor}
\begin{proof}
Let $T$ be jointly strongly stable. It follows from Lemma \ref{uman500} that $\mathfrak{d}_{T}=\mathfrak{d}_{T^*}$. By virtue of \cite[Theorem 5.2, Chapter 6]{SN},  the contraction $T$ belongs to the class $C_0$ defined in Section \ref{sec.2.1}. It follows from Proposition \ref{uman400} that $\chi(T)=0$.

Assume that $\chi(T)=0$. Applying the item (i) of Proposition \ref{uman400} we obtain $\mathfrak{d}_{T}=\mathfrak{d}_{T^*}$. Since $\mathfrak{d}_{T^*}<\infty$, the last equality means that $T$ belongs to~$C_0$~\cite{SN} and, therefore, $T\in{C_{00}}$.
\end{proof} 

The proof of the corollary establishes a somewhat 
stronger result than that stated. 
In fact, we show that $T$ belongs 
to an important subclass $C_0$ of~$C_{00}$. 
The contractions in $C_0$ enjoy a number of interesting properties and admit a convenient characterization. 
Before stating the result, we recall that the characteristic function $\Theta_T(\lambda)$ 
has the \emph{scalar multiply} $\delta(\lambda)$ if $\delta(\lambda)$ is a~scalar-valued analytic function, 
$\delta(\lambda)\not\equiv{0}$, and there exists a contractive analytic function 
$\Omega(\lambda) \colon \mathfrak{D}_{T^*} \to \mathfrak{D}_T$ such that
\begin{equation}\label{uman999}
\Theta_T(\lambda)\Omega(\lambda)=\delta(\lambda)I_{\mathfrak{D}_{T^*}}, \quad
\Omega(\lambda)\Theta_T(\lambda)=\delta(\lambda)I_{\mathfrak{D}_{T}}, \quad \lambda\in\mathbb{D}.
\end{equation}
\begin{lem} [{\cite[Theorem 5.1, Chapter VI]{SN}}] \label{uman777} 
The following are equivalent for $T\in{C_{00}}$:
\begin{itemize}
\item[(i)] $T\in{C_0}$;
\item[(ii)] the characteristic function $\Theta_T(\lambda)$ has the scalar multiply;
\item[(iii)] there exists
a nonzero $\delta(\lambda)\in{H^{\infty}}(\mathbb{D})$ such that $\delta(T)=0$.
\end{itemize}
\end{lem}

\subsection{Examples}

We illustrate the above results by considering Carleson frames and weighted shift operators as generators.
\subsubsection{Remark on Carleson frames}
Consider a bounded operator $T$ in $\Hil$ such that 
$Te_k=\lambda_k e_k$, where $\{\lambda_k\}_{k=0}^\infty$ is a bounded sequence of complex scalars and  $\{e_k\}_{k=0}^\infty$ is an orthonormal basis for $\Hil$. 

Let $\gamma\in\Hil$. The following question naturally arises: under what conditions does the sequence $\{T^n\gamma\}_{n=0}^\infty$ form a dynamical frame? A comprehensive answer was given in \cite{second} (see also \cite{51, third, chri2}):   

\begin{thm}\label{t1}
The sequence $\{T^n\gamma\}_{n=0}^\infty$ is a dynamical frame if and only if the following conditions are satisfied:
\begin{itemize}
\item[(i)] $|\lambda _k|<1$; 
\item[(ii)] $|\lambda_k| \to 1$ as $k \to \infty$; 
\item[(iii)] the sequence $\{\lambda_k\}_{k=0}^\infty$ satisfies the Carleson condition 
\begin{equation}\label{uman234}
\inf_n\prod_{k\ne n}\frac{\vert\lambda_n-\lambda_k\vert}{\vert1-\lambda_n\overline{\lambda_k}\vert}>0; 
\end{equation}
\item[(iv)] the vector $\gamma$ has the form 
$\gamma=\sum_{k=0}^\infty m_k\sqrt{1-|\lambda_k|^2}e_k$,
where $0<c\leq|m_k|\leq{C}$.
\end{itemize}
\end{thm} 

Traditionally, the dynamical frames 
$\mathcal{F}_{T,\gamma}=\{T^n\gamma\}_{n=0}^\infty$ constructed in Theorem \ref{t1} are referred to as \emph{Carleson frames}, while a sequence of complex numbers $\{\lambda_k\}_{k=0}^\infty$ satisfying conditions (i)--(iii) of Theorem \ref{t1} is called a~\emph{Carleson sequence} \cite{Rev, chri4, KM}.

Let $\mathcal{F}_{T,\gamma}$ be a Carleson frame. It follows from (i) and (iii) that the generator operator $T$ must be a contraction in $\Hil$ and that the sequence $\{\lambda_k\}_{k=0}^\infty$ consists of pairwise distinct numbers. 
In the next section, we investigate a~more general setting in which $T$ is assumed to satisfy only these two conditions, 
without necessarily being the generator of a Carleson frame $\mathcal{F}_{T, \gamma}$.

\subsubsection{The contraction $T$.}

Let $\{e_k\}_{k=0}^\infty$  be an orthonormal basis for $\Hil$. Consider a contraction defined by the formula
\begin{equation}\label{uman987}
Te_k=\lambda_k e_k, \quad \lambda_k\in\mathbb{D},  \quad \lambda_k\not=\lambda_j \text{ for } k\not= j, \text{ where } k, j \in\mathbb{N}_0.
\end{equation}
It is easily verified that $T\in{C_{00}}$. Moreover,
\[
(I-TT^*)e_k=(I-T^*T)e_k=(1-|\lambda_k|^2)e_k.
\]
Therefore, $D_{T^*}=D_{T}$, and
\begin{equation*}
    D_Tf=D_{T^*}f=(I-TT^{*})^{\frac{1}{2}}f=\sum_{k=0}^{\infty}\sqrt{1-\vert\lambda_k\vert^2}\langle f,e_k\rangle e_k, \quad f\in\Hil.
\end{equation*}
 This means that $\mathfrak{D}_{T}=\mathfrak{D}_{T^*}=\Hil$ and $\mathfrak{d}_{T}=\mathfrak{d}_{T^*}=\infty$.
The obtained relations show that every Carleson frame $\mathcal{F}_{T,\gamma}$ is defect-space generated. However, it cannot be a Parseval frame, since Theorem~\ref{uuu1212} implies that the set of generating vectors must be infinite.

\begin{prop}\label{uman843} The operator $T$ defined by \eqref{uman987} belongs to $C_0$ if and only if 
the following Blaschke condition holds:
\begin{equation}\label{uman555}
\sum_{k=0}^\infty(1-|\lambda_k|)<\infty.
\end{equation}
\end{prop}
\begin{proof}
If \eqref{uman555} holds, then
the Blaschke product
\begin{equation}\label{uman078}
\delta(\lambda)=\lambda^m\prod_{\lambda_k\not=0}^\infty\frac{|\lambda_k|}{\lambda_k}\frac{\lambda_k-\lambda}{1-\lambda\overline{\lambda_k}}, \quad m\in\{0,1\}.
\end{equation}
belongs to $H^\infty(\mathbb{D})$ and is an inner function. Furthermore, it is easy to see that $\delta(T)e_k=0$, for all $k\in\mathbb{N}_0$. Therefore, $\delta(T)=0$. By Lemma \ref{uman777}, we conclude that $T\in{C_0}$.

Assume that $T\in{C_0}$. Then, according to Lemma \ref{uman777}, the characteristic function $\Theta_T(\lambda)$ has a scalar multiplier $\delta(\lambda)$.  It follows from \eqref{uman987} and Lemma~\ref{uman002} that $\ker\Theta_T(\lambda_k)\not=\{0\}$.  Therefore, in view of \eqref{uman999}, $\{\lambda_k\}_{k=0}^\infty$ is a zero sequence of $\delta(\lambda)$. This yields \eqref{uman555}.
\end{proof}

By Proposition~\ref{uman843}, the assumption $T\in{C_0}$ implies that $|\lambda_n|\to{1}$ as~$n \to \infty$, $\|T\|=1$, and ensures that the vector\footnote{Without loss of generality, we consider the simplest case of generating vectors for Carleson frames $\mathcal{F}_{T, \gamma}$, assuming that $m_k=1$ for all $k\in\mathbb{N}$.}
\begin{equation}\label{bbb3b}
\gamma=\sum_{n=0}^\infty \sqrt{1-|\lambda_n|^2}e_n
\end{equation}
belongs to $\Hil$. On the other hand,  operators $T$ satisfying $\sup_{n\in\mathbb{N}_0}|\lambda_n|<1$ provide examples of operators belonging to the class $C_{00} \setminus C_0$.

\subsubsection{Bessel sequences.}
Assume that $T\in{C_0}$. Another natural question arises: what additional condition must $T$ satisfy for the sequence $\{T^n\gamma\}_{n=0}^\infty$ to be a Bessel sequence in $\Hil$?
Proposition \ref{uman333} provides a simple answer: the sequence $\{T^n\gamma\}_{n=0}^\infty$ is a Bessel sequence if and only if there exists a~nonnegative solution $S$ of the Stein equation
\begin{equation}\label{uman098}
(S-TST^*)f=\langle{f, \gamma}\rangle\gamma, \quad f\in\Hil.
\end{equation}

Assume that $S$ is a solution of \eqref{uman098}. Then, for $f=e_k$, in view of \eqref{uman987} and \eqref{bbb3b},
\begin{equation*}
(S-TST^*)e_k=(I-\overline{\lambda}_{k}T)Se_k=\sqrt{1-|\lambda_k|^2}\sum_{n=0}^\infty\sqrt{1-|\lambda_n|^2}e_n. 
\end{equation*}
Since $(I-\overline{\lambda}_kT)^{-1}e_n=\frac{1}{1-\overline{\lambda}_k\lambda_n}e_n$, we obtain
\begin{equation*}
\begin{split}
Se_k&=(I-\overline{\lambda_k}T)^{-1}\sum_{n=0}^\infty\sqrt{(1-|\lambda_k|^2)(1-|\lambda_n|^2)}e_n\\
&=\sum_{n=0}^\infty\frac{\sqrt{(1-|\lambda_k|^2)(1-|\lambda_n|^2)}}{1-\overline{\lambda}_k\lambda_n}e_n.
\end{split}
\end{equation*}
Thus, the action of $S$ on $\Hil$ is completely determined by the infinite-dimensional matrix 
\begin{equation}\label{uman012}
G=[s_{nk}]_{n,k=0}^{\infty}, \quad s_{nk}=\frac{\sqrt{(1-|\lambda_k|^2)(1-|\lambda_n|^2)}}{1-\overline{\lambda}_k{\lambda_n}}.
\end{equation}
The matrix $G$ is the Gram matrix of the sequence of normalized reproducing kernels in $H^2(\mathbb{D})$; see, e.g., \cite{AAAA, FKP}. In particular,
\begin{equation}\label{uman014}
s_{nk}=\langle{e}_{\lambda_k}, e_{\lambda_n} \rangle_{H^2}, \quad  e_{\lambda_k}(\lambda)=\frac{\sqrt{1-\vert\lambda_k\vert^2}}{1-\lambda\overline{\lambda_k}}, \quad n, k\in\mathbb{N}_0.
\end{equation}

According to \cite[Lemma 14.16]{FKP}, the Blaschke condition \eqref{uman555} indicates that the sequence $\{e_{\lambda_k}\}_{k=0}^\infty$ is minimal and complete in the Hilbert space $K_\delta=H^2(\mathbb{D})\ominus\delta(\lambda)H^2(\mathbb{D})$, where the Blaschke product $\delta(\lambda)$ is defined in $\eqref{uman078}$.
The minimality of $\{e_{\lambda_k}\}_{k=0}^\infty$ (i.e., $e_{\lambda_m}\not\in\overline{\operatorname{span}\{e_{\lambda_k}\}}_{k\not=m}$) implies that every finite Gram matrix $G_N=[s_{nk}]_{n,k=0}^N$ is positive definite. Therefore, the solution $S$ of the Stein equation \eqref{uman098} considered on $\operatorname{span}\{e_k\}_{k=0}^\infty$ is positive.
However, as the following example demonstrates, the Blaschke condition \eqref{uman555} alone does not guaranty the boundedness of $S$.
\begin{ex}\label{ex}
The sequence $\{\lambda_k\}_{k=2}^{\infty}$, where $\lambda_k=1-\frac{1}{k^2}$ for $k\geq2$, satisfies the Blaschke condition. Put $a_k=1/k^2$. Since $1-\lambda_k^2 \geq a_k$ and $1-\overline{\lambda_k}\lambda_n\leq a_k+a_n$
we obtain
$s_{nk}\geq\frac{\sqrt{a_ka_n}}{a_k+a_n}$.
For $N\geq2$ and $N\leq n,k\leq2N-1$, we have
$a_n,a_k\in\left[\frac{1}{4N^2},\frac{1}{N^2}\right]$,
and therefore $s_{nk}\geq\frac18$. Consequently, for
$x_N=\frac1{\sqrt N}\sum_{j=N}^{2N-1}e_j$,
which is a unit vector, we obtain
\[
\|Sx_N\|\geq\langle Sx_N,x_N\rangle
=\frac1N\sum_{n,k=N}^{2N-1}s_{nk}
\geq\frac N8.
\]
Thus the operator $S$ is unbounded.
\end{ex}
The operator $S$ is bounded if and only if the matrix $G$ defines a bounded operator on $\ell^2(\mathbb{N}_0)$. This is equivalent to the sequence $\{e_{\lambda_k}\}_{k=0}^\infty$ being a Bessel sequence in $H^2(\mathbb{D})$; see, e.g., \cite[p. 191]{heil}. When this
sequence is Bessel, its Gram matrix $G$ coincides with $\theta\theta^*$, where 
$\theta \colon H^2(\mathbb{D})\to\ell_2(\mathbb{N}_0)$  is the analysis operator of $\{e_{\lambda_k}\}_{k=0}^\infty$. This means that $G=\theta\theta^*$ is bounded if and only if $\theta$ is bounded. A direct calculation gives 
\begin{equation*}
\theta f=\{\langle f,e_{\lambda_k}\rangle_{H^2}\}_{k=0}^\infty=\left\{\sqrt{1-\vert\lambda_k\vert^2}f(\lambda _k)\right\}_{k=0}^\infty, \quad f\in{H}^2(\mathbb{D}).
\end{equation*}

Define the finite discrete measure $\mu$ on the unit disk by 
\[\mu=\sum_{k=0}^\infty(1-\vert\lambda_k\vert^2)\delta_{\lambda_k},\] where $\delta_{\lambda_k}$ denotes a Dirac delta function. In other words, for any Borel subset $A\subset\mathbb{D}$,
\begin{equation*}
\mu(A)=\sum_{\lambda_k\in{A}}(1-|\lambda_k|^2).
\end{equation*}
The boundedness of $\theta$ then means that there exists $B>0$ such that
\begin{equation}\label{uman014bb}
\|\theta{f}\|_{\ell_2}^2=\sum_{k=0}^\infty(1-|\lambda_k|^2)|f(\lambda_k)|^2=\int_{\mathbb{D}}|f(\lambda)|^2d\mu\leq{B}\|f\|_{{H}^2}^2, \quad f\in{{H^2(\mathbb{D})}}.
\end{equation}
Therefore, $\mu$ is a Carleson measure. Following \cite[p. 614]{FKP} we say that the sequence $\{\lambda_k\}_{k=0}^\infty$ satisfies \emph{Carleson--Newman condition} if \eqref{uman014bb} holds.  This condition is equivalent to the following  \cite[p. 196]{FKP}: 
\begin{equation}\label{uman001}
\sup_{I\subset\mathbb T}
\frac{1}{|I|}
\sum_{\lambda_k\in S(I)}
(1-|\lambda_k|^2)<\infty.
\end{equation}
Here $I=(e^{i\theta_1}, e^{i\theta_2})$, where $0<\theta_2-\theta_1\leq{2\pi}$, is any arc of $\mathbb{T}$ of normalized length $|I|=(\theta_2-\theta_1)/2\pi$, and
$S(I)=\{re^{i\theta} :  1-|I|<r<1, \, \theta_1<\theta<\theta_2\}$
is the corresponding Carleson box.
 
In summary, the operator $S$ is bounded whenever condition \eqref{uman001} holds. Applying Proposition \ref{uman333}  we obtain the following result.

\begin{prop}\label{prop-uman012}
Let $T$ be a contraction defined by \eqref{uman987} and let $\{\lambda_k\}_{k=0}^\infty$ satisfy the  Blaschke condition \eqref{uman555}. Then the following  are equivalent:
\begin{itemize}
\item[(i)] $\{\lambda_k\}_{k=0}^\infty$ satisfies the Carleson--Newman condition \eqref{uman001}; 
\item[(ii)]
the set $\{T^n\gamma\}_{n=0}^\infty$, where $\gamma$ is defined by \eqref{bbb3b}, is a Bessel sequence in~$\Hil$.
\end{itemize}
\end{prop}

\begin{rem}
Another characterization of the Bessel sequences $\{T^n\gamma\}_{n=0}^\infty$ is given in \cite[Lemma 2.6]{48}. 
\end{rem}

\subsubsection{Carleson frame $\mathcal{F}_{T,\gamma}$ and the characteristic function $\Theta_T(\lambda)$.}
The preceding considerations, together with the results of \cite[Chapter 15]{FKP}, allow us to characterize the case in which the Bessel sequence  $\{T^n\gamma\}_{n=0}^\infty$ is a dynamical frame in $\Hil$.

\begin{thm}\label{the1} 
Let a contraction $T$ be defined by 
\eqref{uman987} and let $\{\lambda_k\}_{k=0}^\infty$ satisfy the  Blaschke condition \eqref{uman555}.  Then the  following are equivalent:
\begin{itemize}
\item[(i)]   $\{T^n\gamma\}_{n=0}^\infty$ is a  Carleson frame $\mathcal{F}_{T,\gamma}$ with the generating vector $\gamma$ defined by \eqref{bbb3b};
\item[(ii)]  $\{T^n\gamma\}_{n=0}^\infty$ is a Bessel sequence in $\Hil$ and  
\begin{equation}\label{uman013}
\inf_{k\ne{n}}\|\Theta_T(\lambda_n)e_k\|>0,
\end{equation}
where  $\Theta_T(\lambda)$ is the characteristic function of $T$;
\item[(iii)] the following condition holds:
\begin{equation}\label{uman51}
\inf_n\prod_{k\ne n}\Vert\Theta_T(\lambda_n)e_k\Vert>0.
\end{equation}
\end{itemize}
\end{thm}

\begin{proof}
It follows from Propositions \ref{uman88} and \ref{uman333} that the sequence $\{T^n\gamma\}_{n=0}^{\infty}$ is a dynamical frame if and only if the solution $S$ of the Stein equation \eqref{uman098} is uniformly positive. 
This is equivalent to the fact that the matrix $G$ (see \eqref{uman012}) defines a boundedly invertible operator on $\ell_2(\mathbb{N}_0)$. Equivalently, the sequence of normalized reproducing kernels $\{e_{\lambda_k}\}_{k=0}^\infty$ (see \eqref{uman014}) must be a~Riesz basis for the Hilbert space 
$K_\delta=H^2(\mathbb{D})\ominus\delta(\lambda)H^2(\mathbb{D})$, where $\delta(\lambda)$ is the Blaschke product defined in \eqref{uman078}. 
In view of the classical characterization (see, e.g., \cite[Theorem 15.7]{FKP}), the Riesz basis property of $\{e_{\lambda_k}\}_{k=0}^\infty$ in~$K_\delta$ is equivalent to the Carleson condition 
\eqref{uman234} for the sequence $\{\lambda_k\}_{k=0}^\infty$.
Using~\eqref{uman987}, we~obtain
\[
\begin{split}
\Theta_T(\lambda)e_k&=-Te_k+\lambda D_{T^*}(I-\lambda T^*)^{-1}D_Te_k\\
&=-\lambda_ke_k+
\lambda(1-|\lambda_k|^2)
\sum_{n=0}^{\infty}(\lambda\overline{\lambda_k})^n e_k\\
&=\left(-\lambda_k+
\lambda\frac{1-|\lambda_k|^2}{1-\lambda\overline{\lambda_k}}\right)e_k
=\frac{\lambda-\lambda_k}{1-\lambda\overline{\lambda_k}}e_k.
\end{split}
\]
Thus, the characteristic function of $T$ is given by
\begin{equation*}
\Theta_T(\lambda)e_k=\frac{\lambda-\lambda_k}
{1-\lambda\overline{\lambda_k}}e_k, \quad k\in\mathbb{N}_0, \quad \lambda\in\mathbb{D}
\end{equation*}
and therefore, 
\begin{equation}\label{uman045}
\|\Theta_T(\lambda_n)e_k\|=\frac{|\lambda_n-\lambda_k|}
{|1-\lambda_n\overline{\lambda_k}|}.
\end{equation}
This shows that the Carleson condition \eqref{uman234} can be rewritten as \eqref{uman51}. The~equivalence between (i) and (iii) is established. 

Assume that (ii) holds. In view of Proposition \ref{prop-uman012}, the sequence $\{\lambda_k\}_{k=0}^\infty$ satisfies the Carleson--Newman condition. Furthermore, by means of \eqref{uman045}, 
the condition \eqref{uman013} takes the form $\inf_{k\ne n}\frac{|\lambda_n-\lambda_k|}
{|1-\lambda_n\overline{\lambda_k}|}>0$. This means that the sequence $\{\lambda_k\}_{k=0}^\infty$ is separated \cite[p. 614]{FKP}. Therefore, $\{\lambda_k\}_{k=0}^\infty$  is separated and satisfies the Carleson–Newman condition. By \cite[Theorem 15.3]{FKP}, these two conditions are equivalent to the Carleson condition \eqref{uman234} for $\{\lambda_k\}_{k=0}^\infty$, or equivalently to its reformulation \eqref{uman51}.  The equivalence between (ii) and~(iii) is proved.
\end{proof}

\begin{cor}\label{the2} 
Let $\mathcal{F}_{T, \gamma}$ be a Carleson frame with the generator operator $T$ defined by \eqref{uman987} and the Carleson sequence $\{\lambda_k\}_{k=0}^\infty$. Let 
$\{\mu_k\}_{k=0}^\infty$ be a~sequence of pairwise distinct numbers in~$\mathbb{D}$ satisfying 
\begin{equation*}
\sup_{k\in\mathbb{N}_0}\|\Theta_T(\mu_k)e_k\|\leq\alpha,
\end{equation*}
where $\alpha<1$ and 
\begin{equation*}
\frac{2\alpha}{1+\alpha^2}<\inf_n\prod_{k\ne n}\Vert\Theta_T(\lambda_n)e_k\Vert.
\end{equation*}
Then $\{\mu_k\}_{k=0}^\infty$ is a Carleson sequence, and the family $\mathcal{F}_{T_\mu,\gamma_\mu}=\{T_\mu^n\gamma_\mu\}_{n=0}^\infty$
is a Carleson frame, where 
\[
T_\mu e_k=\mu_ke_k, \quad \gamma_{\mu}=\sum_{k=0}^\infty \sqrt{1-|\mu_k|^2}e_k.
\]
\end{cor}

\begin{proof}
This follows from Theorem \ref{the1}, \cite[Theorem 15.6]{FKP} and the relation~\eqref{uman045}.
\end{proof}

\subsubsection{The case of weighted unilateral shift}
Consider a weighted unilateral shift $T$ defined on vectors of an orthonormal basis $\{e_k\}_{k=0}^\infty$ of $\Hil$:
\begin{equation*}
Te_k=\omega_k e_{k+1}, \quad 0<\omega_k\leq{1},  \quad  k\in\mathbb{N}_0.
\end{equation*}
The condition $0<\omega_k\leq{1}$ guaranties that $T$ is contractive and that its adjoint 
\begin{equation*}
T^*e_0=0, \quad T^*e_{k+1}=\omega_ke_k, \quad k\in\mathbb{N}_0
\end{equation*}
is strongly stable. Moreover, $(I-T^*T)e_k=(1-\omega_k^2)e_k$ and
\begin{equation*}
(I-TT^*)e_0=e_0, \quad (I-TT^*)e_{k+1}=(1-\omega_k^2)e_{k+1}, \quad k\in\mathbb{N}_0.
\end{equation*}
Therefore,
\begin{equation}\label{uman003}
D_Te_k=\sqrt{1-\omega_k^2}e_k, \quad D_{T^*}e_{k+1}=\sqrt{1-\omega_k^2}e_{k+1}, \quad D_{T^*}e_{0}=e_0, \quad k\in\mathbb{N}_0.
\end{equation}

The defect indices $\mathfrak{d}_{T}$ and $\mathfrak{d}_{T^*}$ of $T$ depend on the choice of weights~$\omega_k$.  In the case where $\omega_k=1$ for all $k\in\mathbb{N}_0$ (the classical unilateral shift),  $D_T=0$. Consequently,  $\mathfrak{D}_{T}=\overline{\mathcal{R}(D_T)}=\{0\}$ and $\mathfrak{d}_{T}=0$. On the other hand, the defect operator $D_{T^*}$ coincides with the orthogonal projection of $\Hil$ onto $\ker{T^*}=\operatorname{span}\{e_0\}$ and $\mathfrak{D}_{T^*}=\overline{\mathcal{R}(D_{T^*})}=\ker{T^*}$. 
The last relation means $\mathfrak{d}_{T^*}=\dim\ker{T^*}=1$. In view of Proposition \ref{uman140}, $\gamma(T)=1$. 
The generating vector $\gamma$ can be chosen as $\gamma=e_0$, in which case the dynamical frame $\mathcal{F}_{T,\gamma}$
coincides with the orthonormal basis $\{e_k\}_{k=0}^\infty$. 

Consider a generalization of the preceding case, assuming that
\begin{equation}\label{end14}
\prod_{k=0}^\infty\omega_k=\alpha>0.
\end{equation}
Then
$T^ke_0=\beta_k{e_k}$, where $\beta_k=\prod_{j=0}^{k-1}\omega_k\geq\alpha$ for all $k\in\mathbb{N}$. Therefore,
the~dynamical frame $\mathcal{F}_{T,e_0}=\{T^ke_0=\beta_ke_k\}_{k=0}^\infty$ is a Riesz basis of $\Hil$ with the frame bounds $\alpha^2$ and $1$ and the frame index $\gamma(T)=1$.  The numbers 
\begin{equation*}
\omega_k=\frac{\sqrt{(k+1)(k+3)}}{k+2}, \quad k\in\mathbb{N}_0
\end{equation*}
provide an example of weights satisfying \eqref{end14}. In this case, a direct calculation yields $\alpha=\frac{1}{\sqrt{2}}$ and $D_{T}e_k=\frac{1}{k+2}e_k$, $D_{T^*}e_k=\frac{1}{k+1}e_k$. Consequently, $\mathfrak{d}_{T}=\mathfrak{d}_{T^*}=\infty$.

What happens when $\prod_{k=0}^\infty\omega_k=0$? In this case, the vector $e_0$ cannot serve as a generating vector of $\mathcal{F}_{T, \gamma}$, because $\lim_{k\to\infty}\beta_k=0$ and $\{T^ke_0=\beta_ke_k\}_{k=0}^\infty$ is shown to be a Schauder basis. 
Assume that
$\sup_{k\in\mathbb{N}_0}\omega_k<1$. Then $\|T\|<1$ and, in view of Proposition \ref{uman141}, $\gamma(T)=\mathfrak{d}_{T^*}=\infty$. This means that the sequence 
$\{T^ke_0=\beta_ke_k\}_{k=0}^\infty$ cannot be represented as a dynamical frame $\mathcal{F}_{T, \gamma}$ generated by $T$ using only finitely many generating vectors $\gamma_j$. By Proposition \ref{uman14},  the weighted unilateral shift $T$ can serve as a generator of dynamical frames $\mathcal{F}_{T, \gamma}$ with infinitely many generating vectors $\gamma_j$. The corresponding frames can be constructed using \eqref{umc} and \eqref{uman003}. In particular, choosing $\tau_j=e_{j-1}$, $j\in\mathbb{N}$, in \eqref{umc}, we obtain the generating vectors $\gamma_1=e_0, \gamma_j=\sqrt{1-\omega_{j-1}^2}e_{j-1}$, $j\geq{2}$. \\

\noindent \textbf{Author contributions}
S. Gwizdek contributed to the investigation, formal analysis, and editing of the manuscript. S. Kużel contributed to the conceptualization, investigation, methodology, and preparation of the original draft. W. Majdak contributed to the investigation, formal analysis, and preparation of the original draft. All authors read and approved the final version of the manuscript.

\noindent \textbf{Funding} The work was partially supported by the Faculty of Applied Mathematics, AGH University of Krakow, through statutory tasks funded by the subsidy of the Ministry of Science and Higher Education (Poland).

\noindent \textbf{Conflict of interest} The authors declare that there is no conflict of interest regarding the publication of this paper.

\noindent \textbf{Data availability} No datasets were generated or analyzed during the current study.

\end{document}